\documentclass[10pt]{amsart}
\usepackage[pagebackref]{hyperref}
\usepackage{sagetex}
\usepackage{mathrsfs}
\usepackage{graphicx}
 \usepackage[all]{xy}
  \usepackage{amssymb}	
  \usepackage{color}
  \usepackage{longtable}
\newtheorem*{thm*}{Theorem}
 \newtheorem*{prop*}{Proposition}
\newtheorem*{rmk*}{\textit{Remark}}
\newtheorem*{corr*}{Corollary}
\newtheorem*{lemm*}{Lemma}
\newtheorem{theorem}{Theorem} 
 \newtheorem{lemm}[theorem]{Lemma}
 \newtheorem{prop}[theorem]{Proposition}
\newtheorem{corr}[theorem]{Corollary}

\theoremstyle{remark}

 \newtheorem{exmple}[theorem]{Example}
 \newtheorem{exmples}[theorem]{Examples}
 \newtheorem{rmk}[theorem]{\textit{Remark}}
\numberwithin{table}{section}
\def\lset{\{}  
	\def\rset{ \}}  
	\def\set#1{\lset#1\rset} 
	\def\sett#1#2{\lset #1 \mid  #2 \rset}  
 \def\half{\frac 12}

\def\la{\langle}
\def\ra{\rangle}

\DeclareMathOperator{\Tr}{Tr}
\newcommand{\bC}{{\mathbb{C}}}

\newcommand{\bP}{{\mathbb{P}}}
\newcommand{\bQ}{{\mathbb{Q}}}

\newcommand{\bZ}{{\mathbb{Z}}}
 \usepackage[cal=boondoxo]{mathalfa}

	\newcommand\cO{{\mathcal O}}

\newcommand\germ[1]{{\mathfrak{#1}}}

	\newcommand\gj{{\germ j}}

  \DeclareFontEncoding{LS1}{}{}
    \DeclareFontSubstitution{LS1}{stix2}{m}{n}
              \DeclareSymbolFont{symbols2}{LS1}{stixfrak} {m} {n}
		\DeclareMathSymbol{\operp}{\mathbin}{symbols2}{"A8}
  		
\def\ba{{\mathbf a}}
\def\bx{{\mathbf x}}
  \def\bb{{\mathbf b}}
   \def\ii{{\mathbf i}}
   \def\b1{{\mathbf 1}}
   
\def\mapright#1{\mathop{\vbox{\ialign{
                                ##\crcr
    ${\scriptstyle\hfil\;\;#1\;\;\hfil}$\crcr
 \noalign{\kern2pt\nointerlineskip}
    \rightarrowfill\crcr}}\;}}
\def\into{\hookrightarrow}
\def\onto{\twoheadrightarrow}

\def\chow#1#2{{ \mathsf {Ch}}^{#1}{(#2)}}
\def\chm#1#2{\mathsf{ChM}^{#1}(#2)}
\def\cors#1#2{{\mathsf{Corr}}^{#1}(#2)}
\def\mc#1{\mu_{#1}}

\def\comp{\raise1pt\hbox{{$\scriptscriptstyle\circ$}}}
\def\id{\text{\rm id}} 
 \def\im{\operatorname{Im}}
 \def\Tr{{}^{\mathsf{T}}\kern-0.9pt} 
 \newcommand\note[1]{{\color{blue}#1}} 

\begin{document}
\title{Motivic aspects of  a remarkable class of Calabi--Yau threefolds}
 \author{Gregory Pearlstein \& Chris Peters}
\date{\today}
\maketitle

\section{Introduction}

A  Calabi--Yau manifold is a compact  complex manifold having a trivial canonical bundle. Smooth hypersurfaces  of degree $n+2$ in projective $n$-space give the first examples. Smooth hypersurfaces of $\bP^n$ of different degrees are not.
This changes if one allows certain singularities and defines an
$n$-dimensional  Calabi--Yau variety to be a complex projective $V$-variety 
whose canonical sheaf is trivial.
Recall that a $V$-variety is a complex variety
 admitting a local cover by charts $(U/G_U)$,
$U\subset \bC^n$ (classically) open and $G_U$ a finite group of holomorphic automorphisms  of $U$. 
The kind of examples we are interested in are hypersurfaces in weighted projective 4-space 
$\bP(a_0,a_1,a_2,a_3,a_4)$ that are quasi-smooth and have trivial canonical sheaf.
These notions are explained in \S~\ref{ssec:weights}.
It turns out that there are many Calabi--Yau threefolds of this kind as demonstrated in Sections~\ref{sec:exmple} and \ref{sec:Egypt}.

Our interests in these examples was motivated after considering  two of the families of elliptic 
surfaces  considered in  \cite{PP24}, namely those of even degree  $2c=14$ in $\bP(1,2,3,7 )$ and  of  even
degree $ 22$ in $\bP(1,2,7,11 )$. Both have an equation of the form $H(x_0,x_1,x_2)-x_3^2=0$
and adding a new weight $1$ variable $s$ the resulting threefold $s^{2c}+H(x_0,x_1,x_2)-x_3^2=0$
is a Calabi--Yau threefold $X_{2c}$ admitting a cyclic group of biholomorphic automorphisms
of order $2c$ generated by the automorphism $g$ resulting by multiplying the variable $s$ with a primitive $2c$-th root of unity.
 
Calculations using SAGE revealed that the eigenspace of $g^{c}$  of  the induces action on $H^3(X_{2c},\bQ)$ up to a Tate twist has the same Hodge \note{numbers}   as $H^1(C,\bQ)(-1)$,
where $C$ is the curve on $X_{2c}$ cut out by the codimension 2 subspace $s=x_3=0$.
This turned out not to be a coincidence. Indeed, the quotient 
of  $Y_2=X_{2c}/ \langle g^2\rangle $ is a Fano threefold and \note{we show that} $H^3(Y_2,\bQ)\simeq H^1(C,\bQ)(-1)$.
This leads to a direct proof of the generalized Hodge conjecture for
the Hodge structure on this space.
Moreover, this space is isomorphic to an isotypical component of $H^3(X_{2c},\bQ)$
under the $\bZ/2c\bZ$-action and thus can be interpreted motivically. For the first of these examples this is detailed in Section~\ref{sec:exmple}.

This phenomenon occurs more generally for symmetric Calabi--Yau threefolds $X$  of type $(2c,[A,1,a,b,c])$, $A|2c$, i.e., those given by a weighted polynomials  of the form  
$$ 
 F:= s^{m} +H(x_0,x_1,x_2)- x_3^2,\quad m=2c/A \text{ even},\quad \deg F= 2c.
$$
On such a variety the cyclic group of order $m$ generated by the automorphism $g$ resulting by multiplying the variable $s$ with a primitive $m$-th root of unity. For each divisor $d\not= m$ of 
$m$  the group generated by  $g^{d}$  acts on $X$ with a Fano threefold $Y_d$ as its quotient.
Besides the two previous examples there are many such symmetric Calabi--Yau threefolds, for instance the threefolds of  Tables~\ref{tab:First},\ref{tab:2c-1}, \ref{tab:2c-2} and \ref{tab:misc-ex}.

The   main results we prove in Section~\ref{sec:Main}  are as follows:
\begin{thm*} Let $X$ be symmetric Calabi-Yau  of type $(2c,[A,1,a,b,c])$. Then
\begin{itemize}
\item  Then there is an orthogonal splitting of rational Hodge structures \footnote{$\Psi_{1,m} H^3(X)=0$ since the subspace of $g$-invariants is zero.} 
 $$H^3(X,\bQ)= \Psi_{m, m} H^3(X,\bQ)  \operp \operp^{d\not=m}  \Psi_{d ,m} H^3(X,\bQ) .
  $$
     The first summand contains the transcendental subspace    $H^3(X)_{\rm tr}\subset H^3(X,\bQ)$ and if $Y_d= X/\la g^{m/d} \ra$,
     then  $H^3(Y_d) \simeq \oplus_{e | d }  \Psi_{e,m} H^3(X)$ if $d\ge 2$.
     In particular, $H^3(Y_{2})\simeq  \Psi_{2,m}H^3(X)$.
    \item The $GCH(1,3)$-conjecture   holds for the summands   $\Psi_{d, m} H^3(X,\bQ)$, $d\not=m$; moreover, for $d=2$ 
 the Abel--Jacobi map  $ J(C) \to J(Y_2) =J( H^3(X,\bQ)  ^{g^2})$    is an isogeny.
\item  $X$ admits self-dual Chow-K\"unneth  decomposition and
 the group-action of the cyclic group generated by the action of $g$   on $X$
 induces a further   decomposition 
  \[
  \chm 3 X=( X,   \Psi_{m, m}) \oplus_{d\not=m,2 } (X, \Psi_{d ,m})\oplus (X,\Psi_{2 ,m }). 
    \]
  The first summand contains  the transcendental motive of $X$, the last  summand is isomorphic to $\chm 1 C(-1)$ with third   Chow group $J(C)(-1)$. The other summands $(X, \Psi_{d ,m})$, $d\not=m,2$ are isomorphic to $\chm 3 {Y_d}$.
\end{itemize}
\end{thm*}

 The composition of this note is as follows. In \S~\ref{ssec:weights} background on weighted hypersurfaces is given, in \S~\ref{ssec:GCH-Fano} we review the proof of the generalized Hodge conjecture for the middle cohomology of a Fano threefold, and in  \S~\ref{ssec:chowmots} background on Chow motives is given. The main results are stated and proven in 
 Section~\ref{sec:Main} while in Section~\ref{sec:exmple}  explicit calculations are performed for an  illustrative example of a  threefold of type $(14,[1,1,2,3,7])$. In Section~\ref{sec:Egypt} the types of those
 symmetric Calabi--Yau threefolds are determined for which the weights all divide the degree (so that these are all  deformations of Fermat-type hypersurfaces) and the decomposition announced in the main theorems is explicitly tabulated in Tables~\ref{tab:First}, using SAGE.

Appendix~\ref{app:1}  contains tables giving all possible sums of 5 Egyptian fractions summing up to 1. The first of these  tables is used to find all possible symmetric Calabi--Yau threefolds of Fermat type while Appendix~\ref{app:2} contains the SAGE code we used for computing the occurring types of representations in the middle cohomology of the symmetric Calabi--Yau threefolds in Table~\ref{tab:First}.

 \section{Prepatory material}
 
 \subsection{Weighted hypersurfaces} 
 \label{ssec:weights}

In this subsection we recall some results from the
 literature on hypersurfaces in weighted projective spaces, e.g.
 \cite{DolWPV,IFWings,steen}. Recall that $\bP:=\bP(a_0,\dots,a_n)$ is
 the quotient of $\bC^{n+1}\setminus\set{0}$ under the $\bC^*$-action
 given by $\lambda (x_0,\dots,x_n)=
 (\lambda^{a_0}x_0,\dots,\lambda^{a_n}x_n)$.  One  may   assume
 that $a_0\le a_1\le \cdots\le a_n$.  The affine piece $x_k\not=0$ is
 the quotient of $\bC^n$ with coordinates
 $(z_0,\dots,\widehat{z_k},\dots,z_n)$ by the action of $\bZ/a_k\bZ$
 given on the coordinate  $z_i =x_i/x_k^{(a_i/a_k)}$  by $\rho^{a_i}z_i$,
 where $\rho$ is a primitive $a_k$-th root of unity. Observe that in
 case $a_0=1$, the coordinates $z_j=x_j/x_0, j=1,\dots,n$ are actual
 coordinates on the affine set $x_0\not=0$; there is no need to divide
 by a finite group action.

 Where  the weighted projective space  $\bP$ has   singularities at most along the $k$-codimensional "simplices"
    $L_{j_1,\dots,j_k}=\set{x_{j_1}=\cdots=x_{j_k}=0 }$.  Depending on the weights such a simplex  does have cyclic quotient singularities
 transversal to  it, namely    in case  the set of
 weights that result after discarding $a_{j_1},\dots, a_{j_k}$ are not
 co-prime, say with gcd equal to $h_{j_1,\dots,j_k}$, and then the transversal singularity  type
is 
 \[
 \frac 1{h_{j_1,\dots,j_k}} (a_{0},\dots, \widehat{a_{j_1}}, \dots,
   \widehat{a_{j_k}},\dots,  a_ {n}).
   \] 
This means that these singularities  are the image of
 $0\times \bC^\ell \subset \bC^{k}\times \bC^\ell $, where $\bZ/h\bZ$
 acts on $ \bC^{k}$ by $\rho_h (x_1,\dots,x_k)=
 (\rho_h^{b_1}x_1,\dots,\rho_h^{b_k}x_k)$, $\rho_h$ a primitive $h$-th
 root of unity.   In particular, the vertices are always singular,  and if all weights are pairwise co-prime, these
 are the only singularities. Less stringently, if  any  $n$-tuple from the collection $\set{a_0,\dots,a_n}$ of
weights is co-prime, the only possible singularities occur in
codimension $\ge 2$. We call such weights \textbf{\emph{well formed}}
and in what follows we shall assume that this is the case.

 If $X$ is a degree $d$ hypersurface in weighted projective space
$\bP(a_0,\dots,a_n) $ its \textbf{\emph{type}} is the symbol
$(d,[a_0,\dots,a_n])$. Following \cite{IFWings}, the integer $\alpha_X=d -(a_0+\dots+a_n)$ is called the
\textbf{\emph{amplitude}} of $X$. If the
 corresponding variety $F=0$ in $\bC^{n+1}$ is only singular at the
 origin, the variety $X$  is called  \textbf{\emph{quasi-smooth}}.
This implies that the possible singularities of quasi-smooth
 hypersurfaces come from the singularities of $\bP$. Such a
 hypersurface has at most cyclic quotient singularities, i.e. it is a
 $V$-variety.  To test if $F=0$ is quasi-smooth one uses the Jacobian criterion: the 
 only solution to $\nabla F (\bx)=0$ is $\bx=(x_0,\dots,x_n)=0$.

 The type $(d,[a_0,\dots,a_n])$ of a  hypersurface of degree $d$ in $\bP(a_0,\dots,a_n)$ is called
 \textbf{\emph{well formed}} if the weights $ a_0,\dots,a_n$ are well formed and if
 moreover $h_{ij}=\gcd(a_i,a_j)$ divides $d$ for $0\le i<j\le n$. All
 our examples  are hypersurfaces with  well formed type. In particular, such  hypersurfaces  have at most singularities in codimension $2$,
 and,  moreover,  the divisorial sheaf  $\cO(\alpha_X)$ is precisely the
 canonical sheaf $\omega_X$.

   \begin{exmples}
   We give two examples of quasi-smooth hypersurfaces having  type  $(d,[a_0,\dots,a_n])$.
 \\
 (1) In case all weights $a_j$ divide $d$, the Fermat-type hypersurfaces $\sum x_j^{d/a_j}=0$ are quasi-smooth. It also follows that
 the general hypersurface of such  a type  $(d,[a_0,\dots,a_n])$ is quasi-smooth.
  \\
 (2). Assume that all but one weight, say $a_j$, divide $d$ and that $d= k a_j+ a_\ell$ for some weight $a_\ell, \ell\not=j$, then
 $\sum_{i\not=j} x_i^{d/a_i}+ x_j^k x_\ell=0$ is quasi-smooth. Again,  
 the general hypersurface such type   is quasi-smooth.
   \end{exmples}

 Some Hodge-theoretic results from \cite{steen} are used below, more specifically, as in the non-weighted case the Griffiths' residue calculus
 can be used to find  the Hodge decomposition for the 
quasi-smooth hypersurface $X=V(F)$  in weighted projective space
$\bP(a_0,\dots,a_n)$  in terms of the Jacobian ring $R_F=\bC[x_0,\dots,x_n]/\gj_F$, where
$\gj_F$ is the Jacobian ideal  of $F$.  In particular, with $\Omega_{n}  = \sum (-1)^j x_j dx_0\wedge \cdots dx_{j-1}\wedge\widehat{dx_j}\wedge \cdots \wedge dx_n$, one has
\begin{eqnarray}
\label{eqn:holform}
H^{n,0}(X) &= &\displaystyle \text{Res}\left ( R_F^{\alpha_X} \cdot \frac {\Omega_{n}}  F\right)   ,\\
H^{n-1,1} (X)&= &\displaystyle\text{Res}\left ( R_F^{\alpha_X+\deg F} \cdot \frac {\Omega_{n}}  {F^2} \right).  \label{eqn:seclevel}
\end{eqnarray} 
From this one sees 
\begin{equation}
 h^{n,0}(X) =\dim H^0(X,\omega_X)= \dim H^0(X,\cO(\alpha(X))), \label{eqn:pg}
\end{equation} 
since it is assumed that the symbol of $X$ is well formed.
 \par
As to deformations, we shall use   the following result.
   
\begin{lemm}[\protect{\cite[\S1]{Tu}}] The subspace $\mathsf{Def}_{ {\rm  proj}}$ of the Kuranishi
space of deformations of $X$  within $\bP(a_0,\dots,a_n)$ is smooth with tangent space
canonically isomorphic to $R_F^{d}$. \label{lem:kuranishi}
The Kuranishi family restricted to $\mathsf{Def}_{ {\rm  proj}}$ 
is called the \textbf{\emph{Kuranishi family of type $([d],(a_0,\dots,a_n))$}}.
\end{lemm}

\subsection{The generalized Hodge conjecture for Fano varieties}
\label{ssec:GCH-Fano}

Since the cohomology of a quasi-smooth subvariety  $X$ of weighted projective space has a pure Hodge structure it makes sense 
to consider the generalized Hodge conjecture \allowbreak $GHC(k,n,X)$ for those.  Recall that it states that for every rational Hodge substructure $H'\subset H^n(X)\cap F^kH^n(X)$   one can find a subvariety $Z\subset X$ of codimension $\ge k$ on which $H'$ is supported, i.e., 
$H'=f_* H^{n-2k}(\tilde Z)(-k))$, where $\tilde Z$ is a resolution of singularities of $Z$ and $f_*$ is the Gysin map associated to the natural map $f:\tilde Z\to X$. 
Usually, one calls the subspace of $H^n(X)$ generated by those substructures $ H'$ for fixed $k$ the \textbf{\emph{subspace $N^kH^n(X)$ of co-level $k$}}, while 
the largest  rational Hodge substructure $N^k_{\rm{Hdg}}H^n(X)\subset H^n(X)\cap F^kH^n(X)$ is called the   \textbf{\emph{subspace of Hodge level $k$}}.  So $GHC(k,n,X)$ states $N^kH^n(X)=N^k_{\rm{Hdg}}H^n(X)$.

In our situation we shall be considering   $GHC(1,3,X)$ for $X$ a threefold. In this case the conjecture is equivalent to the existence of a smooth projective
surface $S$ (not necessarily irreducible) admitting a morphism $f: S\to X$ such that $f_* H^1(S,\bQ)= N^1_{\rm{Hdg}}H^3(X,\bQ))$.
Note also that $N^1_{\rm{Hdg}}H^3(X,\bQ))$ is the smallest rational Hodge substructure of $H^{2,1}(X)\oplus H^{1,2}(X)$.
For a Fano threefold $X$  by definition minus the canonical divisor is  ample. So $H^{3,0}(X)=0$ and then $N^1_{\rm{Hdg}}H^3(X,\bQ))=H^3(X,\bQ)$. This also holds if $X$ is a quasi-smooth hypersurface of weighted projective 4-space whose type is well formed and with negative amplitude.
Such $X$ is an example of a so-called \textbf{\emph{$\bQ$-Fano threefold}}.
The validity of $GHC(1,3,X)$ in this case is well known, but for completeness we sketch the simple proof.

\begin{prop} \label{prop:ghc} Let $X$ be $\bQ$-Fano threefold. Then $GHC(1,3,X)$ holds.
\end{prop}
\begin{proof}[Sketch of the proof.] The crucial ingredient here is that $X$ is uniruled, that is, $X$ is covered by rational curves, as  for example 
shown in \cite{RCVars}. Using this, one follows the strategy of Conte and Murre in \cite{conte}. The 3 steps of their proof apply in this case:

\begin{itemize}
\item the conjecture  $GHC(1,3,X)$ is stable under morphisms of finite degree,
\item it is stable under birational morphisms,
\item it holds for $\bP^1\times S$, where $S$ is a surface.
\end{itemize}
 See  the proof of \cite[Proposition 13.3]{lewis} in Lewis's monograph for the first two points. The third is obvious since
 $N^1_{\rm{Hdg}}H^3(\bP^1\times S,\bQ))=H^3(\bP^1\times S,\bQ))$ and $N^1_{\rm{Hdg}}H^3(X,\bQ))= i_* H^1(S,\bQ)(-1)$  
 \end{proof}

 \subsection{Chowmotives}
 \label{ssec:chowmots}
\begin{small} 
In this subsection the basics of Chow motives  is recalled as explained more fully in \cite{andre2,motifbk}. In fact,  it is slightly extended to the category of projective $V$-varieties.
 \end{small}
 
 \medskip 
The categorical nature of motives comes from  correspondences between varieties $X$ and $Y$, that is, cycles on their product.
For Chow motives one considers these up to rational equivalence, in other words in the Chow group of $X\times Y$.
In fact, one needs more, namely a product structure on the Chow groups which it  into a ring. This is classically possible
for smooth projective varieties, but here it is used for  projective  $V$-variety such as a quasi-smooth subvariety of a weighted projective space.
This forces one to pass to the Chow groups with $\bQ$-coefficients as explained in Mumford's basic article~\cite{enumgeom}. 
Consequently, for $X$   a  projective  $V$-variety the notation 
$$\chow{} {X}=\oplus_{r=0}^ {\dim X}\chow r X
$$ 
stands for  the Chow ring with $\bQ$-coefficients.
Moreover, for a correspondence $\Gamma$ from  a  $V$-variety  $X$ to a $V$-variety $Y$, its class $[\Gamma]\in \chow{} {X\times Y}$
 shall also be called a \textbf{\emph{correspondence}}.  It  is said to have degree $r$ it it belongs to
\[
\cors r {X,Y}:=  \chow {r+\dim X} {X\times Y}.
\]
For example,  the graph  of  a morphism $f:X\to Y$ defines a  correspondence 
$\Gamma_f \in \cors r {X,Y}$ where the degree equals $r=\dim Y-\dim X$, while its transpose  in $ \cors 0{Y,X}$
has degree $0$.  If $\Gamma \in \cors r {X,Y}$ there are induced homomorphisms
\[
\Gamma_* : \chow i X \to \chow{i+r} Y,\quad \Gamma_* : H^i(X,\bQ) \to H^{i+2r}(Y,\bQ),
\]
 in particular, if $\Gamma=\Gamma_f$, the induced operations coincide with the usual homomorphisms $f_*$ on Chow groups and cohomology,
 while the action of its transpose corresponds to the usual induced homomorphisms $f^*$.

Correspondences can be composed. In particular, a self-correspondence  
 $ p \in \chow{} {X\times X}$ is called a \textbf{\emph{projector}}, if $p\comp p=p$. For instance, the diagonal $\Delta_X$ of $X$ is a projector.
 A (pure) \textbf{\emph{Chow motive}} $(X,p)$ consists of an  projective  $V$-variety together with a projector $p$. 
Projectors have degree $0$ and morphism $(X,p) \to  (Y,q)$  between Chow motives   by definition belong to  $q\comp \cors 0  {X,Y} \comp p$. 
Chow motives admit direct sums and images and kernels.  For instance if $M=(X,p )$, then a projector $q$ of $M$ is an element $q=p\comp q'\comp p$ with $q'\in \cors 0  {X,X}$ such that $q\comp q=q$ and $N=(X,q)$ is the image of $q$. Note that $q=p\comp q =q\comp p$, i.e., $N$ is a
 constituent  of $M$. So is $N'=\ker (q)=(X,p-q)$, and we have $M=N\oplus N'$. 
 
Instead of pure motives, degree $m$ motives are triples $(X,p,m)$ where $(X,p)$ is a pure motive and $m\in\bZ$.
The degree is only used to change the notion of a morphism   $(X,p,m) \to  (Y,q,n)$: it is an element of $q\comp \cors {n-m}  {X,Y} \comp p$.
Motives admit a \textbf{\emph{tensor product}}
\[
 (X,p,m)\otimes (Y,q,n) := (X\times Y,p\times q, m+n),
\]
with $ \b1= \chm {} {\text{pt}}$.
Degree $m$ motives can always be obtained from pure motives upon  tensoring with  the weight $m$ Tate motive $\mathbf T^{\otimes m}$,
$\mathbf T= (\text{pt},\id,1)  $:
\[
(X,p,m)\simeq (X,p,0)(m):=  (X,p,0)\otimes \mathbf T ^{\otimes m}.
\]
 Motives have their Chow groups and cohomology groups:
 \[
 \chow i {X,p}:= \im \left(p_*: \chow i X \to \chow i X\right),\quad H^i(X,p):= \im ( p_* ) \subset H^{i }(X) .
 \]

The \textbf{\emph{Chow motive of a   projective    V-variety $X$}}   by definition  is the pair
$\chm{} X := ( X,\Delta_X)$.
It  has 2 natural constituents,  defined by two projectors
\[
p_0(X):= x \times X,\quad p_{2d}(X)= X\times x, \quad x\in X,\, d=\dim X,
\]
that is,
\[
\chm 0  X:= (X,p_0(X)),\quad \chm {2d} X = (X,p_{2d}(X)).
\]
The two projectors are orthogonal in the sense that  $ p_{2d}\comp p_0=p_0\comp  p_{2d}=0$ and then $p^+(X)=\Delta_X- p_0-p_{2d}$
is also a projector and one has a direct sum decomposition
\[
\chm{}X=  \chm 0  X \oplus \chm +   X\oplus \chm {2d}   X,\quad \chm +   X=(X, p^+(X).
\]
As an example, the \textbf{\emph{Lefschetz motive}} is defined as
$\mathbf L= \chm 2 {\bP^1}$ which under $\otimes$ is the inverse of the Tate motive, that is $\mathbf L\otimes \mathbf T= \b1 $ .

As to morphisms of motives, note that the usual morphisms $f:X\to Y$ leading in general not to degree $0$ correspondences,
the motivic morphism $\chm {}f :\chm{} Y\to \chm{} X$ associated to $f$ is the  degree $0$ correspondence   given by the transpose of $\Gamma_f$.  
  
We say that a $V$-variety $X$ admits a \textbf{\emph{Chow-K\"unneth decomposition}} (C-K
decomposition for short) if there exist orthogonal projectors $p_i(X)\in  \cors{0} {X,X}$ for $0\le i \le 2d$  decomposing the
diagonal of $X$, i.e., 
\[
p_i(X)\comp p_j (X)= \left\{ \begin{array}{cc}
0 & j\ne i \\
p_i (X)& j=i
\end{array}
\right.,\quad  \sum_{i=0}^{2d} p_i (X)= \Delta_X,
\]
 and such that, moreover the cohomology class $p_i(X)^*$ belongs to  the K\"unneth component $ [\Delta_X]^i\in H^{2d-i}\otimes H^i(X)$
of the cohomology class of the diagonal. If  the projectors
can be chosen such that
 $
p_{2d-i} (X)= {\Tr p}_i(X)$  the C-K decomposition is said to be 
\textbf{\emph{self-dual}}.  This uses the  \textbf{\emph{dual $(X,p,m)^*$ of a motive}}, given by $(X,p,m)^*:= (X, \Tr p, d-m)$.
Hence the Chow motive of $X$ decomposes as 
\[
\chm{}X= \oplus_{j=0}^{2d} \chm j X,\quad 
 \chm j X (d) \simeq (\chm {2d-j} X)^* .
 \]

 \begin{exmples} \label{exm:CD}
 (1). Surfaces $S$ admit a self-dual Chow-K\"unneth decomposition (see e.g. \cite[\S 6.3]{motifbk}). Moreover, $\chm 2 S$
splits into an algebraic motive $\mathsf A(S)$  isomorphic to a direct sum of Lefschetz motives and a transcendental motive $\mathsf T(S)$.
These are characterized by their cohomology groups: the first has cohomology the subgroup spanned by the algebraic classes and the
cohomology group of the second  consists of the transcendental cycles.
 \\
 (2). 
 Suppose that $X   $ is an projective $V$-variety  of dimension $d$ for  which  $H^i(X)$ is algebraic for all $i\ne d$, i.e., $H^{2j}(X)$, $2j\not=d$
 is generated by classes of algebraic subvarieties and $H^i(X)=0$ for $i\not=d$ odd. 
 Then $X$ admits a self-dual Chow-K\"unneth decomposition.  
 For a proof consult \cite[Appendix D]{motifbk}.
  \end{exmples}
 
 \begin{rmk} \label{rmk:transmot}
We will be considering quasi-smooth weighted threefold hypersurfaces.
 For those, the Hodge decomposition is $H^3(X,\bC)= H^{3,0}(X)\oplus H^{2,1}(X)\oplus H^{1,2}(X)\oplus H^{0,3}(X)$.
 The smallest rational sub-Hodgestructure $T (X)\subset H^3(X,\bQ)$ whose complexification contains $H^{3,0}(X)$ is called
 the \textbf\emph{transcendental cohomology}, while its orthogonal complement (under cup product) is a rational Hodge structure $A(X)$ which
 contains the subspace  $N^1H^3(X) \subset H^3(X,\bQ)$ supported on a surface, i.e.  $N^1H^3(X)= \bigcup _Z i_* H^1(Z,\bQ)$ where
    $Z $ is a smooth surface with a generic embedding   $i : Z \to X$. The generalized Hodge conjecture would imply that  $A(X)=N^1H^3(X)$.
   Even if we don't know this, suppose that  $\chm 3 X= \mathsf T(X) \oplus \mathsf A(X)$, where $H^3(\mathsf T(X))= T(X)$ and
   $H^3( \mathsf A(X))=A(X)$ we shall call $\mathsf T(X)$ the \textbf{\emph{transcendental motive of $X$}}.
 \end{rmk}

\subsection{Group representations of finite cyclic groups and motives}
\label{ssec:grpreps}

Let $\mc m$ be the cyclic group of order $m$ with generator $ g$. For each divisor $d$ of $m$ there is exactly one irreducible
representation of degree $\phi(d)$, the degree of the   $d$-th  cyclotomic polynomial.  
These irreducible representations decompose the  group ring $\bQ[\mu_m]$. For each divisor $d$ of $m$ this determines a projection  $\Psi_{d,m} :\bQ[\mu_m]\to \bQ[\mu_m]$ which is the identity on the given representation and $0$ on all all other direct summands.
This then gives a decomposition of the identity into orthogonal idempotents
\begin{eqnarray}
\Psi_{d,m} \comp \Psi_{d',m}  =  \begin{cases}
						0 & \text{ if } d\not= d'\\
						\Psi_{d,m} & \text{ if } d = d'.
\end{cases} \label{eqn:OrthIdemp}
\quad & \sum_{d|m} \Psi_{d,m}  =   1.  
\end{eqnarray} 
We shall identify these projectors with the corresponding representations, which actually are their images in the group ring. As projectors they can be given as an element of the group ring considered as an operator, i.e.,

\begin{equation}
\label{eqn:IrrRep}\
\Psi_{d,m} =\sum_{k=0}^{m-1}  a_k g^k\in \bQ[\mu_m].
\end{equation} 
To find the decomposition of any rational $\mc m$-representation one may apply the following result. 
\begin{lemm}\label{lemm:isotyp}
A   $\mc m$-representation $\bQ$-vector space $V$ splits into  a direct sum of its isotypical components $\Psi_{d,m}V$.\\
1. If $V_k$ is the eigenspace of the action of a generator $g$ of  $\mc m$
on $V$ with eigenvalue $\rho_m^k$, one   has 
$$
\Psi_{d,m} (V\otimes \bC) := \oplus_{k, \gcd(k,d)=1} V_{k(m/d)}.
$$
2. The subspace of $V$ on which $g^d$ acts trivially is equal to the direct sum
 $ \oplus_{e | d}  \Psi_{e,m} V$.
 
\end{lemm}
\begin{proof}  Recall the expression for the cyclotomic polynomial of degree $d$ for $d$ a divisor of $m$:
\begin{eqnarray*}
\Phi_d(x) &= &\prod_{1\le k\le d,\, (k,d)=1} (x- \rho_d^k)\\
		&= &\prod_{1\le k\le d,\, (k,d)=1} (x- \rho_m ^{k(m/d) })
\end{eqnarray*}
Hence the corresponding direct sum $\Psi_{d,m}V$
of the eigenspaces $V_k\subset V\otimes\bC$
is defined over $\bQ$.  This proves 1.
\\
2. Let $H$ be the group generated by $g^d$. Then
$ V^H= \oplus _{e, e |d } V_{(m/d)\cdot e} =  \oplus_{e | d}  \Psi_{e,m} V$ 
as    follows from the definitions. 
 \end{proof}

\begin{exmples}  
1. 
 If $p$ is a prime, $\mu_p$ one has only two irreducible representations, the trivial one, $\b1$, and the representation $\Psi_{p,p }$ of degree  $p-1$.
Then   \eqref{eqn:IrrRep} reads  
 \begin{eqnarray*} 
 \b1   =  \frac 1 p[ 1+ \sum_{j=1}^{p-1}  g^j ], \qquad 
 \Psi_{p,p }  = \frac  1p [ (p-1)-   \sum_{j=1}^{p-1}  g^j ].
\end{eqnarray*}
 The first  expression for $\b1$ is clear and the second is the unique projector orthogonal to $\b1$ and summing up to the identity. 
\\
2. Consider $\mu_{2^k}$. It has irreducible representations $\b1$ of degree 1 and  $\Psi_{2^\ell,2^k}$, $\ell=1,\dots,k$ of degree $2^{\ell-1}$.  
For instance, for $\mu_8 $ one has 
\begin{eqnarray*}
 \Psi_{2,8 }& =& \frac 18 (1-g+g^2-g^3+g^4-g^5+g^6-g^7 ), \quad\dim  \Psi_{2,8 } =1 \\
\Psi_{4,8 } &= & \frac 18 (2-2g^2+2g^4-2g^6), \quad\dim  \Psi_{4,8 } =2\\
\Psi_{8,8 } &= & \frac 18 (4 -4g^4), \quad\dim  \Psi_{8,8 } =4.
\end{eqnarray*}

\end{exmples}

Suppose $\mu_{m}=\la g\ra $ acts as a group of morphisms on a $V$-variety $X$. Then   the induced action $g^*$ on $H^*(X,\bQ)$ preserves the Hodge structure.
 All eigenspaces  $H^r(X,\bC)_k$    inherit a decomposition
$\oplus H^r(X,\bC)_k^{p,q}= H^r(X,\bC)_k\cap H^{p,q}(X)$, $p+q=r$, but these need not be
a Hodge structures. However the iso-typical components of $H^j(X,\bQ)$, denoted 
$\Psi_{d,m}H^j(X,\bQ)$ 
are indeed Hodge structures since these are  real Hodge structures with underlying
$\bQ$-structure.

Since the $\Psi_{d,m}$, $d|m$ form orthogonal idempotents the action  of $\mu_{m}$ on the diagonal  
  $\Delta_X\subset X \times X$, viewed as a projector,  decomposes into an  orthogonal
sum of  projectors
\[
\Delta=\sum_{d|m} \Delta_{d,m},\quad  \Delta_{d,m}= \Psi_{d,m}\comp \Delta,
\]
and hence  the motive of $X$ decomposes as $\chm {} X =\oplus_{d|m}  (X, \Delta_{d,m} )$.
 
\begin{exmple} \label{exm:ChowCD}
Let $X$ be a quasi-smooth hypersurface  of weighted projective $n+1$-space which is stable under the action of $\mu_m$.
By Example~\ref{exm:CD}, $X$ admits self-dual C-K decomposition  $\chm {}X=\oplus_{j=0}^{2n} (X,p_j)$ which further decomposes under the action of $\mu_m$.This action is non-trivial only on $\chm n X$ and there one has 
 $\chm n X= \oplus_{d|m} (X,\Delta_{d,m} ) $.
\end{exmple}

 \section{Main results}
 \label{sec:Main}
 
 Our main interest concerns  quasi-smooth Calabi--Yau hypersurfaces $X$ in weighted 4-space $\bP(A,1,a,b,c)$, $1\le a\le b< c$,  of degree $ 2c$.
 Since we assume $(2c,[A,1,a,b,c])$  to be well formed, $\omega_X=\cO (2c-(A+1+a+b+c))$. This sheaf being trivial, one has  $A= c- (1+a+b)\ge 1$.
 Moreover, we demand that \textbf{\emph{$c$ is divisible by $A$}}.
 As we shall see in Section~\ref{sec:Egypt}, there are many such threefolds.

 Our results do not concern     all of them, but rather those whose equation
 in homogeneous coordinates\footnote{This is chosen in order to attune the notation to that used in \cite{PP24}.}$(s, x_0,x_1,x_2,x_3)$  has the following specific form
 \begin{equation}
 \label{eqn:CYs}
 F:= s^{m} +H(x_0,x_1,x_2)- x_3^2,\quad m=2c/A,\quad \deg F= 2c.
 \end{equation} 
 Such   a Calabi--Yau threefold  in   will be called \textbf{\emph{symmetric Calabi--Yau of type $(2c,[A,1,a,b,c])$}}.
 Indeed,  such threefolds have an action by the cyclic group $\mc  m$: the generator $g$ sends $(s:x_0:x_1;x_2:x_3)$ to
 $(\rho_{m}s:x_0:x_1:x_2:x_3)$ with $ \rho_{m}$ a primitive $m$-th root of unity. 
 
 \begin{rmk} \label{rmk:OnModuli} The projective moduli of the symmetric Calabi--Yau threefold   of type $(2c,[A,1,a,b,c])$ is  one more  
 than the projective moduli of  the family of curves with type  $(2c,[1,a,b])$,
 since the variables $s$ and $x_3$ can be scaled so that these come with coefficient $1$, which still allows  to scale $H$.
 \end{rmk}

 The quotients of $X$ by the action of the proper subgroups  of  $\mc {m}$ of order $d$ generated by $g^{m/d}$  all
 are  $\bQ$-Fano varieties    
   \begin{equation}
 \label{eqn:modell}
 Y_d= V(G_d)\subset \bP(A m/d,1,a,b,c),\, G_d:= t^{d} +H(x_0,x_1,x_2)- x_3^2,\quad t=s^{m/d}. 
 \end{equation} 
 of type $(2c, [Ad,1,a,b,c])$.
 If $X$ is quasi-smooth, then so is $Y_d$.  Since \allowbreak $(2c,[Ad,1,a,b,c])$  is also well formed,  $\omega_{Y_d}=\cO( A(1-d) )$, and so $Y_d$ is a $\bQ$-Fano threefold.
 
 The  threefold  $Y_{2}$ plays a special role: the action of $g$ on $X$ induces an involution on $Y_{2}$ whose quotient is $\bP(1,a,b,c)$ with branchlocus the surface $S=V(H)$ given by
 the equation $H(x_0,x_1,x_2)+ x_3^2=0$. Setting $x_3=0$ produces the curve $C=V(H)\subset \bP(1,a,b)$ whose symbol $(2c,[1,a,b])$ is well-formed and so   the genus $g(C)$ can be calculated from equation~\eqref{eqn:pg}. In fact, using Griffiths' residue calculus as given in equations \eqref{eqn:holform} and \eqref{eqn:seclevel}, we have

 \begin{lemm} $g(C)= h^{2,1}(Y_{2})$.
 \end{lemm}
 \begin{proof} Since $\gj_{Y_{2}}$ is generated by $t, x_2, \gj_C$, the Jacobian  ring
 of $Y_{2}\in \bP(c, 1,a,b,c)$ is isomorphic to that of $C$. In particular, the parts of degree  $2c-(1+a+b)$ agree.
 For $Y_{2}$ the dimension equals $h^{2,1}(Y_{2})$, while for $C$ its dimension equals the genus,
 since the symbol $(2c,[1,a, b])$ of $C$ is and so $\omega_C=\cO(2c-(1+a+b)$.
 \end{proof}

  In the diagram below the relation between the various varieties is depicted.
 
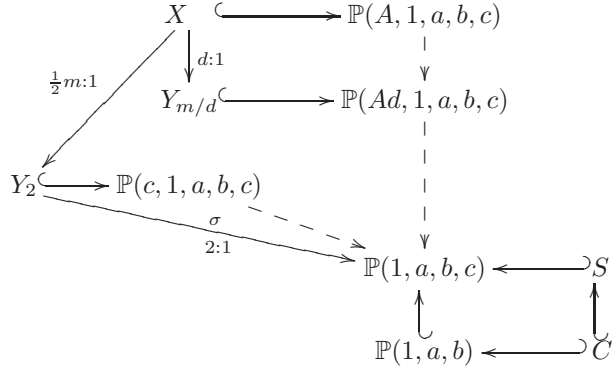
\xymatrix@R=1.3pc{\quad  & X \quad \ar[ddl]_{\half m :1}\ar@{^{(}->}[r]  \ar[d]^{d :1}    & \bP(A,1,a,b,c)  \ar@{-->}[d]\\
&  Y_{m/d}  \ar@{^{(}->}[r] &   \bP( Ad ,1,a,b,c) \ar@{-->}[dd]  \\
Y_{2}   \ar@<-0.5ex>[drr]_{2:1}^{\sigma\quad} \ar@{^{(}->}[r] &   \bP(  c,1,a,b,c) \ar@{-->}[dr]   \\
                          &                           & \bP( 1,a,b,c)                 &   S \ar@{_{(}->}[l] \\
                          &  &\bP(1,a,b)\ar@<0.5ex>@{_{(}->}[u]     &   \, C\ar@<0.5ex>@{_{(}->}[u]  \ar@{_{(}->}[l]
}
\noindent For each point 
 $\ba=(0:a_0:a_1:,a_2:a_3:0)\in C \subset \bP(c,1,a,b,c)$ the line
 \begin{equation}
  L_\ba=\sett{(\lambda:\mu \cdot \ba:\lambda)}{(\lambda:\mu)\in \bP^1} \label{eqn:L}
 \end{equation} 
 belongs to $Y_{2}$ and the union gives the surface 
$T_{Y_{2}}= \bigcup_{\ba\in C}  L_\ba  \subset Y_{2} $ 
which is a cone on $C$ with vertex $(1:0:0:0:1)$.
Now the crucial observation is as follows.

\begin{prop}\label{prop:IJ}
Let $i: C\times \bP^1 \onto T_{Y_{2}}\into Y_{2}$ the natural immersion. Then the Gysin map $i_*: H^1(C\times \bP^1,\bZ)(-1)\to  H^3(Y_{2},\bZ)$ is an isomorphism
of Hodge structures.
\end{prop}
\begin{proof} It suffices to show that the Gysin map is injective since both source and target have dimension $g=g(C)$. 
If $\set{a_1,\dots,a_{g},b_1,\dots, b_{g}}$ is a standard symplectic basis for $H_1(C)$,
 it suffices to prove that the  3-cycles  $A_i, B_j$, $i,j\in \set{1,\dots,g}$
swept out by $L_\ba$ when $\ba$ traverses $a_i,b_j$   give independent
homology classes in $Y_{2}$. Now note that the line   $L^-_\ba\subset Y$ passing through $\ba$ and $(1:0:0:0:-1)$ only meets $L_\ba$ in
$\ba $.  Taking cycles $a'_i$ in $C$ disjoint but homologous to $a_i$ and similarly $b'_j$ disjoint but homologous to $b_j$,
the  3-cycles  $A'_i, B'_j$, $i,j\in \set{1,\dots,g}$ swept out by these by the lines $L^-_\ba$ when $\ba$ traverses the homologous basis $a'_i,b'_j$ for 
$H_1(C,\bZ)$ meet the cycles $A_i,B_j$ transversely showing that $A_i\cdot B'_j=\delta_{ij}$, $A_i\cdot A'_j=0=B_i\cdot B'_j$. Hence 
 if for some rational numbers $r_i,s_j$ there is a relation $\sum r_i [A_i]+ \sum_j  s_j [B_j]=0$
between the classes of these 3-cycles, then intersecting with all $ A'_i $ and $  B'_j $ shows that the relation is trivial, proving that $i_*$ is injective.
\end{proof}

\begin{corr}  \label{corr:OnHalfm} Let $Y_{2}$ be as in \eqref{eqn:modell}. Then 
\begin{enumerate}
\item The generalized Hodge conjecture for
$H^3({Y_{2}})$ holds.
\item The Abel--Jacobi map $ J(C) \to J({Y_{2}})$ is an isogeny.
\end{enumerate}
\end{corr}
\begin{proof}
(1). The lemma implies  that    the entire cohomology $H^3({Y_{2}},\bZ)$ is  carried  by $i_*H^1(T,\bZ)$. This proves the generalized Hodge conjecture for
$H^3({Y_{2}})$.\\
(2) Recall that the Abel--Jacobi map $ \alpha$ is given by
$[\omega] \mapsto \int_\Gamma i_* (p^* \omega)$, where $\omega$ is a holomorphic 1-form on $C$ and 
where $\Gamma$ a 3-cycle whose boundary equals  $ [L_ \ba ] - [L_ \bb]$.
Recall also that the tangent map   of $\alpha$ at $0$ is just $ i_* \comp p^* $. Since this induces a Hodge-isometry of $\bQ$-Hodge structures, the
induced morphism $\alpha$ is an isogeny.
\end{proof}

 The Hodge structure for  $H^3(X,\bQ)$ can be  interpreted  in terms of its structure as a  $\mu_{m}$-representation.
 First consider $H^{3,0}(X,\bC)$ which is 1-dimensional spanned by the residue of $\Omega_{\bP(A,1,a,b,c)}/F$. From \eqref{eqn:holform}
 one sees that this is the eigenspace for the eigenvalue $\rho_{m}$.  By Lemma~\ref{lemm:isotyp} the corresponding $\bQ$-representation space 
 $\Psi_{m,m}H^3(X,\bQ)$ 
 is then found by adding all eigenspaces whose eigenvalues are of the form $\rho^k_{m}$ with $\gcd(k,m)=1$.
 In particular, all other representation spaces $\Psi_{d, m} H^3(X,\bQ)$, $d|m$, $d\not=m$ are either empty 
 or of pure Hodge type $(2,1)+(1,2)$.
 In fact, they are all $\bQ$-Hodge structures (since $g$ acts as a holomorphic automorphism of $X$) and so 
 $h^{2,1}( \Psi_{d, m} H^3(X,\bQ))=h^{1,2}(\Psi_{d, m }H^3(X,\bQ))$. Indeed, since this is the invariant part of  the $\mc {m/d}$-action on $X$ given by $g^{d}$,  one has
 \[
 \Psi_{d, m} H^3(X,\bQ) = H^3(Y_d,\bQ),\quad d\not= m.
 \]
 Applying also part 2 of Lemma~\ref{lemm:isotyp}   this yields one of the main results:

 \begin{theorem} \label{thm:main1} Let $X$ be symmetric Calabi-Yau  of type $(2c,[A,1,a,b,c])$, i.e. $X=V(F)$ with $F$ as in \eqref{eqn:CYs}. The cyclic group $\mu_m$ with   $m=2c/A$ 
 generated by $(s:x_0:x_1:x_2:x_3)\mapsto (\rho_{m}s:x_0:x_1:x_2:x_3)$, $\rho_{m}$ a primitive $m$-th root of unity, acts on $X$.
 Then, with the notation of \S~\ref{ssec:grpreps}  there is an orthogonal splitting of rational Hodge structures 
 $$
 H^3(X,\bQ)= \Psi_{m, m} H^3(X,\bQ)  \operp \operp^{d|m,d\not=   m}  \Psi_{d ,m} H^3(X,\bQ) .
  $$
     The first summand contains the transcendental subspace    $H^3(X)_{\rm tr}\subset H^3(X,\bQ) $ and {$H^3(Y_d) \simeq \oplus_{e | d, e }  \Psi_{e,m} H^3(X)$ if $d\ge 2$}.
     In particular, $H^3(Y_{2})\simeq  \Psi_{2,m}H^3(X)$.
    \end{theorem}

Consider now  the point $\ba =(0,a_0,a_1,a_2,a_3,0)\in C \subset \bP(A,1,a,b,c)$. Then the line
  $L'_\ba=\sett{(\lambda\cdot ( 1:0:0:0:1)+  \mu \cdot \ba }{(\lambda:\mu)\in \bP^1}$ belongs to $X$ which gives  the surface 
  $T_X= \bigcup_{\ba\in C} L'_{\ba}$   in $X$.
The $\mu_{\half m}$-action on $X$ given by $g^2$, where $g:X\to X$ is as in \eqref{eqn:Acts}, fixes the point  $P= (1:0:0:0:1)\in X$ and hence also $T_X$.
  By Proposition~\ref{prop:ghc} and Corollary~\ref{corr:OnHalfm}  one thus has:
 
 \begin{corr} (1) The $GCH(1,3)$-conjecture   holds for the summands   $\Psi_{d, m} H^3(X,\bQ)$, $d\not= m$. \\
 (2) For $d= 2$ the Abel--Jacobi map 
$ J(C)\to J(Y_{2}) =J( H^3(X,\bQ)  ^{g^2})$    is an isogeny.
\end{corr}

\begin{rmk}  It is not clear from this approach whether the generalized Hodge conjecture  is true for the summand that contains the transcendental part   (as it should). 
  \end{rmk}

Using the results from Examples~\ref{exm:CD} and \ref{exm:ChowCD} one can upgrade Theorem~\ref{thm:main1} to a result about Chow motives.

\begin{theorem}\label{thm:main2} Let $X$ be a symmetric C-Y  of type $(2c,[A,1,a,b,c])$. Then 
\\
 (1) $X$ admits self-dual C-K decomposition. \\
 (2) The group-action of $\mc {m}$ on $X$
 induces a further   decomposition 
  \[
  \chm 3 X=( X,   \Delta_{m, m} ) \oplus_{d\not=m, 2} (X,  \Delta_{d ,m})\oplus (X,\ \Delta_{2 ,m }). 
    \]
  The first summand contains  the transcendental motive of $X$, the last  summand is isomorphic to $\chm 1 C(-1)$ with third   Chow group $J(C)(-1)$. 
  Moreover, if $d\not=m $, then 
  $\chm 3 {Y_d} \simeq \oplus_{e| d} (X, \Psi_{e ,m} \comp\Delta)$.

 \end{theorem}

\section{An explicit example}
\label{sec:exmple}

There  are two examples of symmetric Calabi--Yau threefolds constructed from the elliptic surfaces  considered in  \cite{PP24},
namely those of type   $(14,[1,2,3,7])$ and $ (22,[1,2,7,11])$. Both have amplitude 1 and so
$(14,[1,1,2,3,7])$ and $(22,[1,1,2,7,11])$ give symmetric Calabi--Yau threefolds  provided we choose their  equation as in \eqref{eqn:modell}.

 For simplicity we shall only give detailed calculations  for the first case and with  the choice   $H=x_0^{14}+x_1^7-x_2^4x_1$, that is,
$F:=  s^{14}+x_0^{14}+x_1^7-x_2^4x_1-x_3^2$. Observe that the type of $V(H)\subset \bP(1,2,3)$ is NOT well-formed, but since it is
a quasi-smooth curve,  Griffiths' residue calculus resulting in   \eqref{eqn:holform} and \eqref{eqn:seclevel} can be applied and shows that its genus is 10.
Applying this calculus to $V(F)$ one finds also: 
\\
 (1)  $F^3=H^{3,0}$ is 1-dimensional with basis the class of the residue of the rational 4-form $\Omega_4/F$ with pole along $V(F)$. 
 \\ 
(2) $F_2/F_3\simeq H^{2,1}$ has a $132$-dimensional basis,  the class of the residue of the rational 4-forms $M\cdot \Omega_4 /F^2_{14}$,
where $M$ runs over a basis for the degree 14 part of $R/\gj_F$.
The jacobian ideal $\gj_F$ has monomial generators
$s^{13},  x_0^{13},  x_1^6- x_2^4, x_3, x_1x_2^3$ and $x_1^7 $ giving a monomial basis for the degree 14 part of $R/\gj_F$, where $R=\bC[s,x_0,x_1,x_2,x_3]$ as in Table~\ref{tab:MonBasis}. 

\begin{rmk}
By Remark~\ref{rmk:OnModuli} that we may vary $H$ in a family having  $19 $ projective moduli, while    the full family
has  $132=h^{1,2}(V(F))$ projective moduli. A similar result holds for the threefold of type $(22,[1,1,2,7,11])$: one finds $19$, respectively $214$ moduli.
\end{rmk}

The $\mu_7$-action by $g^ 2$ produces a quotient Fano threefold  $Y_2\subset \bP(7,1,2,3,7)$ 
with equation $G_2=t^2+x_0^{14}+x_1^7-x_2^4x_1-x_3^2=0$ (with coordinates 
$t=s^7$, $y_i=x_i, i=0,\dots,3$) and which has middle Hodge numbers $(0,10,10,0)$. Its quotient by the remaining involution $\sigma$  is $\bP(1,2,3,7)$ so that  $\sigma^*$ must act  as $-\id$ on $H^3(Y_2)$. This follows also from the table since $\gj_{G_2}$
is generated by $t,  x_0^{13},  x_1^6- x_2^4, x_3, x_1x_2^3$ so that  a monomial basis for the degree 14 part of $R/\gj_{G_2}$
is also a basis for  the degree 7 part of $(\bC[x_0,x_1,x_2]/\gj_C) $ since this  is also obtained by the monomials obtained by the entries in Table~\ref{tab:MonBasis} corresponding to a monomal $M$ exactly divisible by $s^6$   after dividing by $s^6$.
\par
The involution on $X$ given by $g^{7}$ produces the Fano threefold $Y_7\subset \bP(2,1,2,3,7)$
with equation $ u^7+x_0^{14}+x_1^7-x_2^4x_1-x_3^2=0$, $u=s^2$ and has middle Hodge numbers $(0,60,60,0)$ as indicated by the blue stars in Table~\ref{tab:MonBasis} below.

In Table~\ref{tab:MonBasis}  we have marked the $-1$-eigenspace of $H^{2,1}X )$ of $g^*$ where
\begin{equation}
g :X\to X,\quad g(s:x_0,x_1,x_2,x_3)\mapsto (\rho_{14}\cdot s:x_0,x_1,x_2,x_3),\, \rho_{14}= \exp(2\pi\ii /14). \label{eqn:Acts}
 \end{equation} 
 It corresponds to
the occurrence of $s^6$ and is denoted by $*$.
So this eigenspace has dimension $10$ inside a 132-dimensional subspace.
Dividing the  monomials in the third column of Table~\ref{tab:MonBasis} below for $*=7$  by   $s^6$ we get a monomial basis for $(R'/\gj_C)^{8} \simeq H^{1,0}(C ) $, where 
\begin{equation}
 C \subset \bP(1,2,3),\quad X\cap \set{s=x_3=0}. \label{eqn:C}
\end{equation}

\begin{table}[htp]
\caption{Monomial basis for $H^{2,1}(X)$} \label{tab:MonBasis}

\begin{center}
\begin{tabular}{|l|l|c|}
\hline
monomials $M$ &range  &  $*$ of any color:$s^{*-1}$  occurs in $M $  \\
&&  $*\in ({\color{blue}1},{\color{red} 2},{\color{blue}3},{\color{red} 4},{\color{blue}5},{\color{red} 6},{\color{black}7},{\color{red}8},{\color{blue}9},{\color{red}10},{\color{blue}11},{\color{red}12},{\color{blue}13}) $ \\
\hline
$s^kx_0^{14-k}$ & $k=2,\dots,12$          &    $(0,0,{\color{blue}*},{\color{red}*},{\color{blue}*}, {\color{red}*} ,* ,{\color{red}*},{\color{blue} *},{\color{red}*},{\color{blue}*},{\color{red}*},{\color{blue} *}) $ \\
$s^kx_0^{12-k}x_1$ & $k=0,\dots,12$     &   $({\color{blue}*}  ,{\color{red}*},{\color{blue}*},{\color{red}*},{\color{blue}*},{\color{red}*} ,*,{\color{red}*},{\color{blue}*},{\color{red}*},{\color{blue}*},{\color{red}*},{\color{blue}*}) $ \\
$s^kx_0^{10-k}x_1^2$  &$ k=0,\dots,10$ &   $({\color{blue}*},{\color{red}*},{\color{blue}* },{\color{red}*},{\color{blue}*},{\color{red}*} ,*,{\color{red}*},{\color{blue}*},{\color{red}*},{\color{blue}*},0,0)$  \\
$s^kx_0^{8-k}x_1^3$ & $ k=0,\dots,8$    &   $({\color{blue}*},{\color{red}*},{\color{blue}*},{\color{red}*},{\color{blue}*},{\color{red}*} ,* ,{\color{red}*},{\color{blue}*},0,0,0,0)$  \\
$s^kx_0^{6-k}x_1^4$ & $ k=0,\dots,6$   &    $({\color{blue}*},{\color{red}*}, {\color{blue}*},{\color{red}*},{\color{blue}*},{\color{red}*} ,*,0,0,0,0,0,0)$ \\
$s^kx_0^{4-k}x_1^5$ & $ k=0,\dots,4$   &    $({\color{blue}*},{\color{red}*},{\color{blue}*},{\color{red}*},{\color{blue}*},0,0,0,0,0,0,0,0)$\\
$s^kx_0^{2-k}x_1^6$ & $ k=0, 1,2$   &  $({\color{blue}*},{\color{red}*},{\color{blue}*},0,0,0,0,0,0,0,0,0,0)$ \\
$s^kx_0^{11-k}x_2$ & $ k=0,\dots,11$  &   $({\color{blue}*},{\color{red}*},{\color{blue}*},{\color{red}*},{\color{blue}*},{\color{red}*} ,* ,{\color{red}*},{\color{blue}*},{\color{red}*}{\color{blue},*},{\color{red}*},0) $ \\
$s^kx_0^{8-k}x_2^2$ & $ k=0,\dots,8$  &   $({\color{blue}*},{\color{red}*},{\color{blue}*},{\color{red}*},{\color{blue}*}, {\color{red}*} ,*    ,{\color{red}*},{\color{blue}*},0,0,0,0) $\\
$s^kx_0^{5-k}x_2^3$ & $ k=0,\dots,5$  &    $({\color{blue}*},{\color{red}*},{\color{blue}*},{\color{red}*},{\color{blue}*},{\color{red}*},0,0,0,0,0,0,0) $   \\
$s^kx_0^{9-k}x_2x_1$ & $ k=0,\dots,9$  &   $({\color{blue}*},{\color{red}*},{\color{blue}*},{\color{red}*},{\color{blue}*}, {\color{red}*},* ,{\color{red}*},{\color{blue}*},{\color{red}*},0,0,0)$ \\
$s^kx_0^{7-k}x_2x_1^2$ & $ k=0,\dots,7$  & $({\color{blue}*},{\color{red}*},{\color{blue}*},{\color{red}*},{\color{blue}*},{\color{red}*} ,*  ,{\color{red}*},0,0,0,0,0)$\\
$s^kx_0^{5-k}x_2x_1^3$ & $ k=0,\dots,5$ &  $({\color{blue}*},{\color{red}*},{\color{blue}*},{\color{red}*},{\color{blue}*},{\color{red}*},0,0,0,0,0,0,0)$\\
$s^kx_0^{3-k}x_2x_1^4$ & $ k=0,\dots,3$ &   $({\color{blue}*},{\color{red}*},{\color{blue}*},{\color{red}*},0,0,0,0,0,0,0,0,0)$\\
$s^kx_0^{1-k}x_2x_1^5$ & $ k=0,1$    &   $({\color{blue}*},{\color{red}*},0,0,0,0,0,0,0,0,0,0,0)$\\
$s^kx_0^{6-k} x_2^2 x_1 $ & $ k=0,\dots,6$  &   $({\color{blue}*},{\color{red}*},{\color{blue}*},{\color{red}*},{\color{blue}*},{\color{red}*},*,0,0,0,0,0,0)$\\
$s^kx_0^{4-k} x_2^2 x_1^2 $ & $ k=0,\dots,4$  &   $({\color{blue}*},{\color{red}*},{\color{blue}*},{\color{red}*},{\color{blue}*},0,0,0,0,0,0,0,0)$\\
$s^kx_0^{2-k}x_2^2 x_1^3 $ & $ k=0,1,2 $ &   $({\color{blue}*},{\color{red}*},{\color{blue}*},0,0,0,0,0,0,0,0,0,0)$ \\
$x_2^2x_1^4$&& $({\color{blue}*},0,0,0,0,0,0,0,0,0,0,0,0)$ \\
\hline
\end{tabular}
\end{center}
 \end{table}

The irreducible $\bQ$-represen\-tations of $\mu_{14}$
   are  $\b1=\Psi_{1,14},-\b1=\Psi_{2,14},   \Psi_{7,14 },  \Psi_{14,14}$,   corresponding  to the divisors of $14$.    By lemma  \ref{lemm:isotyp}, the decomposition of $H^3(X,\bQ)$ into isotypical representations can be read off from the eigenspaces
  of $g^*$ on $H^3(X,\bC)$. In this case a basis of $H^{3,0}(X)$ is   the residue of $\Omega_4/F$ and a basis 
  for $H^{2,1}(X)$ is given by the residues of the forms $s^k M_{14-k} (x_0,x_1,x_2)\cdot \Omega_4/F^2$ with eigenvalue  $\rho_{14}$ 
 where $M$ runs through the monomials of Table~\ref{tab:MonBasis}. These all give eigenvectors with eigenvalue $\rho_{14}^{k+1}$.
 To find the eigenvectors on $H^{1,2}(X)\oplus H^{0,3}(X)$ one just has to take the conjugates.
 Then the  eigenspaces for eigenvalues $\rho^k_{14}$, $\gcd(k,14)=1$ give the transcendental eigenspace $ \Psi_{14,14}H^3(X)$.
 Those with $k$ odd give $  \Psi_{7,14}H^3(X)$ and the remaining ones (with eigenvalue $-1$) give  $ \Psi_{2,14}H^3(X)$.
This then gives Table~\ref{tab:MonBasis} where $ \Psi_{14,14 }$  corresponds to the blue stars, $ \Psi_{7,14 }$ to the red stars while the black stars correspond to $-\b1$.  From the table one thus gets the Hodge numbers so that
 the splitting of Theorem~\ref{thm:main2} in terms of their Hodge vectors reads as follows:
 
 \begin{lemm} There is an orthogonal splitting of rational Hodge structures 
 $$H^3(X,\bQ)=  \Psi_{14,14}H^3(X) \operp   \Psi_{7,14}H^3(X) \operp
  \Psi_{2,14}H^3(X)
  $$
   with Hodge vectors  $(1,62,62,1)$, $(0,60,60,0)$, $(0,10,10,0)$ respectively.  The last summand is isometric to $H^3(Y_2,\bQ)$ (since it is the $ \mu_{7}$-invariant part of $H^3(X,\bQ)$).
 \end{lemm}
 
  \section{Fermat type  examples using Egyptian fractions}
\label{sec:Egypt}

In this section a complete list of  symmetric Calabi--Yau threefolds  of Fermat type and of   type $(2c,[A,1,a,b,c])$ are given.
Having this type is equivalent   to $a, b,A$ dividing  $2c$, $A= c-(1+a+b)$ and $2c/A$ being even. 
  Since then  $1 + a + b + c + A = 2c$, dividing by $2c$ gives
  \[
  \frac 1 {2c} + \frac 1 x + \frac 1 y + \frac 12 +\frac 1 t  =1,\quad x=2c/b, \, y= 2c/a,  \, t= 2c/A,
  \]
  an expression of $5$ Egyptian fractions summing up to $1$ with $t$ even. 

 Permuting $(a,b,c)$ so that $a\le b\le c$ this gives  Table~\ref{tab:First} enumerating all 
 101  Fermat-type symmetric weighted Calabi-Yau threefolds
 $X_{a,b,c}$ of degree $2c$ in $\bP(A,1,a,b,c)$ with equations
 $s^{2c/A}+ x_0^{2c}+ x_1^{2c/a}+ x_2^{2c/b}+x_3^2=0$. The   column "$g$" gives the genus of the curve $X_{a,b,c}\cap \set{s=x_3=0}$. The representations $H^3(\bQ)$ are multiples of the irreducible representations of $\mu_m$
 and are tabulated by the divisors $d$ of $m$ and collected by the occurring multiplicities. If $d=m$ the representation
 is the only  sub-Hodge structure of level 3. The contribution is calculated from the dimensions of the eigenspaces 
 for the eigenvalues $\rho_m^k$ with $(k,m)=1$ taking care of this extra contribution. The level 1 types  of irreducible representations 
 can be found in a similar way but one only needs to take care of $(h^{2,1},h^{1,2})$.
 For example, in the second example 12.(48) means 12 copies of an irreducible representation of rank $\phi(48)=16$
 and so this has Hodge numbers (1,95,95,1) while 12.(24,16,12,8,6,4,3,2) has rank
 $12.(\phi(24)+\phi(16)+\phi(12)+\phi(8)+ \phi(6)+\phi( 4)+\phi(3)+\phi(2 ))=372$ and so has Hodge numbers (0,186,186,0).
 Hence  the Hodge numbers  (1,281,281,1)  of the entire middle cohomology,   matching the third column.

 \begin{longtable}{|  l  |  l  |c|c|c|c|c||c| c |  c | l  |}
\caption{Symmetric Calabi-Yau 3-folds of Fermat type} \label{tab:First}
\\
\hline
$(x,y,  t,  2 c)$   & type   & $h^{1,2} $ & $g$& order  & $H^3(\bQ)$ as a representation  \\
\hline 
\endfirsthead%
\caption[]{(continued)}\\
\hline 
$(x,y,  t,  2 c)$   & type   & $h^{1,2} $ & $g$& order  & $H^3(\bQ)$ as a representation  \\
\hline 
\endhead

(3,7,44,924)    & (924,[21,1,132,308,462])& 257& 6& 44& 12.(44,22,11,4,2) \\
(3,7,48,336]    &(336,[7,1,48,112,168])& 281&6 & 48&   12.(48,24,16,12,8,6,4,3,2)\\
(3,7,56,168)   &(168,[3,1,24,56,84]) &329 &6 &56 &12.(56,28,14,8,7,4,2) \\
(3,7,84,84)  &(84,[1,1,12,28,42])&491&6&84& 12.(84,28,21,14,7,4,3,2)\\
&&&&& +11.(42)\\
*(3,25,8,600)& (600,[75,1,24,200,300])&167& 24& 8 &48.(8,4,2)\\
(3,8,26,312)&  (312,[12,1,39,104,156])&174&7&26& 14.(26,13,2) \\
(3,8,28,168)    &(168,[6,1,21,56,84])& 188&7&28&14.(28,14,7,4,2)\\
(3,8,30,120)      &(120,[4,1,15,40,60])&201&7&30&14.(30,15,10,6,5,2)+13.(3)\\
(3,8,32,96) &(96,[3,1,12,32,48])&216&7&32& 14.(32,16,8,4,2)\\
(3,8,36,72)  &(72,[2,1,9,24,36 ])&241&7&36& 14.(36,18,9,3,6,4,2)+13.(12)\\
 (3,8,48,48) &   (48,[1,1,6,16,24])  & 321 &7&   {48}&  
 14.(48,8,6,4,2) \\
 &&&&  &   + 13.(24,12,3) \\
(3,9,20,180)&(180,[9,1,20,60,90])&150&7&20&16.(20,10,5,4)+14.(2)\\
(3,9,24,72)&(72,[3,1,8,24,36])&181&7&24& 16.(24,12,8,4,3)+15.(6)\\
&&&&  & \qquad  +14.(2)\\

  (3,9,36,36) &   (36,[1,1,4,12,18]) &271& 7&   {36}  & 16.(36,12,9,4,3)+\\
  &&&&& \qquad 15.(9)+14.(18,2) \\
   (3,10,16,240)	&(240,[15,1,24,80,120])&134&9&16&  18.(16,8,4,2)\\
*(3,16,10,240)	&(240,[24,1,15,80,120])&134&15&10& 30.(10,5,2)\\

  (3,10,18,90)&(90,[5,1,9,30,45])&151&9&18& 18.(18,9,6,2)+ 17.(3)\\
  *(3,18,10,90)	& (90,[9,1,5,30,45])&151&16&10& 34.(10,5)+32.(2)\\
  (3,10,20,60) &  (60,[3,1,6,20,30]) &  170 & 9&    {20}&   18.(20,10,4,5,2)  \\
 *(3,20,10,60) & (60,[6,1,3,20,30]) &170& 19&   {10}  & 38.(10,5,2) \\
    (3,10,30,30)&   (30,[1,1,3,10,15]) &251& 9&   {30} &  
    17.(30,15,3)+18.(10,5,6,2) \\
     
 (3,12,16,48)  & (48,[3,1,4,16,24]) & 161&10 &    {16}& 22.(16,8)+20.(4,2)\\
    *(3,16,12,48)&  (48,[4,1,3,16,24]) &161& 15&    12   &  29.(12,3)+30(6,4,2)\\
  
 (3,12,12,60)	& (60,[5,1,4,20,30])& 151 &  13& 12 & 28.(12,4,3)+27.(6)\\
 &&&&&\qquad +26.(2)\\
  
(3,12,14,84)	&(84,[6,1,7,28,42])&141&10&14&22.(14,7)+20.(2) \\
(3,12,18,36)& (36,[2,1,3,12,18])&182&10&18&22.(18,9)+21.(6)+20.(3,2)\\

(3,12,24,24)  & (24,[1,1,2,8,12]) &242& 10  &    {24} &  
                  22.(24,8) +20.(12,4,2)\\
                  &&&&& \qquad+21.(6) \\
                    *(3,24,12,14) & (24,[2,1,1,8,12) &242 &22 &  12    &  44.(12,4,3,2)+45.(6)\\  
(3,13,12,156)	&     (156,[13,1,12,52,78])  &131&12&12&24.(12,6,4,3,2)\\
(3,13,24,168)		& (168,[21,1,6,56,84]) &188&27& 8& 54.(8,4,2)\\ 
 (3,14,12,84)&     (84,[7,1,6,28,42]) &141&13&12&26.(12,6,4,2)+25.(3)\\
    (3,14,14,42)&   (42,[3,1,3,14,21]) &168 &13&    {14}& 26.(14,7,2)  \\ 
     
    (3,18,18,18)&  (18,[1,1,1,6,9]) &272 &32 &  {18}  &  32(18,9,3,2)+33.(6) \\
     *(3,18,12,36)  & (36,[3,1,2,12,18]) &182&16 &  12 &  34.(12,4)+33.(6)\\
 &&&&& \qquad    +32.(3,2)\\
 (3,26,8,312)	&(312,[39,1,12,104,156])&174&25&8& 50.(8,4,2)\\
(3,27,8,216)	&(216,[27,1,8,72,108])&180&25&8& 52.(8,4)+50.(2)\\
   (3,30,8,120)		& (120,[15,1,4,40,60]) &201& 28&8&58.(8,4)+56.(2) \\ 
  (3,30,10,30) & (30,[3,1,1,10,15]) &251&28&   {10} &  56.(10,5,2) \\
  (3,32,8,96)      	&(96,[12,1,3,32,48]) &216&31&8&62.(8,4,2) \\ 
  (3,36,8,72)    	& (72,[9,1,2,24,36])  &241&34&8& 70.(8)+68.(4,2) \\ 
 (3,36,8,48)     	& (48,[6,1,1,16,24])   &321&46&8& 92.(8,4,2)\\ 
 (4,5,21,420)	&(420,[105,1,20,84,210])&119&40&4&80.(4,2) \\
  (4,5,24,120)   	& (120,[5,1,24,30,60])  &137&6&24&12.(24,12,8,6,4,3,2) \\ 
   *(5,24,4,120)	&      (120,[30,1,5,24 ,60])&137&46&4& 92.(4,2) \\ 
  (4,5,22, 220)	&(220,[55,1,10, 44,110])&125&42&4& 84.(4,2) \\
   *(5,4,22,220) 	&(220,[10,1,44,55,110])&125&6&22&12.(22,11,2)\\
  (4,5,30,60) & (60,[2,1,12,15,30]) &171 &6 &    {30}&  12.(30,15,10,6,2)+ 11.(5)  \\
  (4,5,40,40) & (40,[1,1,8,10,20]) &227 &6 &     {40}  &  12.(40,10,8,4,2)\\
     &&&&& \qquad\qquad+11.(20,5)\\
(4,6,13,156) 	&(156,[26,1,12,39,78]) &89&18&6&36.(6,3,2) \\
*(6,13,4,156)	&(156,[39,1,12, 26,78])&89&30&4& 60.(4,2)\\
* (4,15,6,60) 	 &(60,[10,1,4,15,30]) &103&21&6&  36.(6,3,2)\\

 (4,6,14,84)&    (84,[6,1,14,21,42]) & 96 & 7&   14 &   15(14,7)+14.(2)\\ 
    *(4,14,6,84)   &    (84,[14,1,6,21,42]) &96& 19&   6 & 39.(6,3)+38.(2) \\
   *(6,14,4,84) &  $((84,[21,1,6,14,42])$ &96 &32 & 4 & 65.(4)+64.(2) \\ 
 (4,6,15,60)& $(60,[15,1,4,10,30])$&103 & 34  & 4 & 70.(4)+68.(2)\\  
(4,6,16,48)&    (48,[3,1,8,12,24])   &110 &7&   16 & 15.(16,8)+14.(4,2)\\
 *(4,16,6,48)&  (48,[8,1,3,12,24])   &110 &21 &    6 &   45.(6,3)+42.(2)\\
 *(16,4,6,28)&  (48,[12,1,3,8,24])   &110 &37&   4 &   74.(4,2)\\
(4,6,18,36)    &    (36,[2,1,6,9,18])   &124 &7 &   18 &  15.(18,9)+14.(6,3,2) \\
 *(4,18,6,36)    &    (36,[6,1,2,9,18])   &124 &25 &     6&  50.(6,3,2) \\
(4,6,24,24)&   (24,[1,1,4,6,12])  &164 &7 &     24 & 15.(24,8)+13.(12)+14(6,3,2)\\
  *(6,24,4,24)   &  (24,[6,1,1,4,12])  &164 &55 &      4  &   110.(4,2)  \\

   *(4,24,6,24)   &  (24,[4,1,1,6,12])  &164 &33 &    6 &   66.(6,3,2) \\
(4,8,12,24)  &(24,[3,1,2, 6,12]) &111& 15 &     8 &  33.(8)+31.(4)+30.(2) \\ 
 
   *(8,12,4,24) &$(24, [6,1,2,3,12])$ & 111& 37 &     4 &     75.(4)+74.(2)  \\
(4,8,16,16) &(16,[1,1,2,4,8])  &147&  9&    16&  21.(16)+19.(4)+ 18.(8,2) \\
   *(4,16,8,16)  &(16,[2,1,1,4,8])  &147& 21  &    8 &   42.(8,2) +43.(4)  \\
 
   *(8,16,4,16)  &(16,[4,1,1,2,8])  &147& 49 &   4   &   99.(4)+98.(2)  \\  
    (4,8,10,40)  &(40,[4,1,5,10,20])  &92& 9&   {10}& 21.(10,5)+18.(2) \\
    *(4,10,8,40)  &(40,[5,1,4,10,20])  &92& 13 &    8&  27.(8)+ 26.(4,2)  \\
    *(8,10,4,40)  &(40,[10,1,4,5,20])  &92& 31 &    4   &   62.(4,2)  \\
(4,10,10,20) &(20,[2,1,2,5,10]) &116& 13 &    {10} & 26.(10,5,2)  \\
   *(10,10,4,20)  &(20,[5,1,2,2,10]) &116& 36 &    4   &   81.(4)+ 72(2) \\
(4,12,12,12)  &  ($12,[1,1,1,3,6])$ &165& 15&   12  &  30.(12,6,3,2)+31.(4) \\
    *(12,12,4,12)           & $(12,[3,1,1,1,6]) $& 165 & 55    &     4  &  111.(4)+110.(2) \\
(4,7,10,140) 	&(140,[14,1,20,35,70])  &80&9&10&18.(10,5,2)\\ 
(4,7,14,28) 	&(28,[2,1,4,7,14])&113&9&14& 18.(14,2)+17.(7)\\ 
(4,8,12,24) 	&(24,[2,1,3,6,12]) &111&9&12& 21.(12,6,3)+19.(4)+18.(2)  \\
(4,8,4,72) 	&(72,[9,1,8,18,36]) &83&12&8& 24.(8,4,2)\\ 
5,5,12,60) 	&(60,[5,1,12,12,30]) &85&6&12&16.(12,6,4,3)+ 12.(2) \\
 (5,5,20,20) & (20,[1,1,4,4,10]) &143& 6&    {20} &  16.(20,4,5)+13.(10) \\
     &&&&& \qquad\quad+12.(2)\\
 (5,6,10,30) & (30,[3,1,5,6,15]) &87& 10&   {10}  &   20.(10,2)+19.(5) \\
     *(5,10,6,30)	&(30,[5,1,3,6,15])&87&16&6& 36.(6,3)+32.(2)\\ 
(5,6,8,120) 	&(120,[15,1,20,24,60]) &69&10&8&20.(8,4,2) \\
  *(5,8,6,120)	&(120,[20,1,15,24,60])&69&14&6&28.(6,3,2)\\ 
 (5,10,10,10) & (10,[1,1,1,2,5]) &145 &16&    {10}   &  33.(10)+ 32.(5,2) \\
 (5,25,4,100)	&(100,[25,1,4,20,50]) &141&46&4& 96.(4)+92.(2) \\
   (5,30,4,60)	&(60,[15,1,2,12,30])&171&56&4& 116.(4)+112.(2)\\ 
 (5,40,4,40)	&(40,[10,1,1,8,20])&227&76&4& 152.(4,2)\\ 
 (6,6,8,24)   & (24,[3,1,4,4,12])  &84& 10  &    8&  25.(8,4)+20.(2)  \\
  *(6,8,6,24)  & (24,[4,1,3,4,12])  &84& 17 &    6 &  34.(6,3,2) \\
(6,6,12,12)   & (12,[1,1,2,2,6]) &126 & 10  &  12 &  25.(12,4)+21.(6,3)+20.(2) \\
   *(6,12,6,12)                 &  ($12,[2,1,1,2,6])$ &126&25&   6& 51.(6,3)+50.(2)  \\
(6,6,7,42) & (42,[7,1,6,7,21]) &74 &15 &   {6}     &  30.(6,3,2)  \\  
(6,6,9,18)  &(18,[3,1,2,3,9]) &95& 19        &   {6}   &  39.(6)+38.(3,2) \\
(6,18,4,36)	&(36,[9,1,2,6,18])&124&40&4& 85.(4)+80.(2)\\
(7,10,4,140)	&(140,[35,1,14,20,70])&80&27&4& 54.(4,2)  \\
(7,14,4,28) & (28,[7,1,2,4,14]) &113& 36&    {4}  &  78.(4)+72.(2)  \\
(8,8,8,8)  &(8,[1,1,1,1,4]) &149& 21&   8&  43.(8,4)+42.(2)\\
 (8,9,16,72)& (72,[18,1,8,9,36) &83&28&   {4}  & 56.(4)+56.(2)  \\  
 (9,9,4,36)& (36,[9,1,4,4,18])&  91& 28&    {4}  & 64.(4)+ 56.(2)  \\ 
    
 \hline
 \end{longtable}%

 \begin{rmk} Threefolds from this list having different types cannot be isomorphic as follows
 from the results of Esser~\cite[Theorem 2.1]{esser}.\\
 \end{rmk}
 
 \section{Some non-Fermat examples}

 In this section, we systematically study the non-Fermat examples of symmetric
symmetric Calabi-Yau hypersurface of degree $2c$ in $\bP(A,1,a,b,d,c)$.  

\subsection*{Case 1:  $a$ and $b$ are divisors of $2c$ and $d$ is a divisor of $2c-1$}
Write $ap=bq=2c$ and $rd=2c-1$.  Then, dividing
$
   1 + a + b + d + c = 2c
$
by $2c$ yields
\begin{equation}
    \frac{1}{2c} + \frac{1}{p} + \frac{1}{q} + \frac{d}{2c} + \frac{1}{2} = 1
    \label{eq:fermat-equality}
\end{equation}
Now, 
 $   
      \displaystyle  \frac{rd}{2c} = 1-\frac{1}{2c}$ which implies $ \displaystyle \frac{d}{2c} = \frac{1}{r} - \frac{1}{2rc}
 $
and so \eqref{eq:fermat-equality} becomes
$$
    \left(\frac{1}{2c}\right)\left(\frac{r-1}{r}\right) + \frac{1}{p} + \frac{1}{q} + \frac{1}{r} + \frac{1}{2} = 1
$$
In particular, since $(r-1)/r < 1$ it follows that we must have 
$$
    1<\frac{1}{2c} + \frac{1}{p} + \frac{1}{q} + \frac{1}{r} + \frac{1}{2}
$$
\begin{table}[h]\caption{$rd=2c-1$}  
\begin{tabular}{|l|l|l|l|l|}
\hline
(1,1,1,5,8) & (2,1,2,3,8) & (1,1,2,7,11) & (2,1,2,9,14) & \underline{(1,1,9,7,18)}  \\
(9,1,1,7,18) & \underline{(3,1,9,5,18)}  & (9,1,3,5,18) & (6,1,4,7,18) & (6,1,6,5,18)  \\
\underline{(1,1,5,13,20)} & (5,1,1,13,20) & \underline{(2,1,4,13,20)}  & (4,1,2,13,20) & \underline{(2,1,14,11,28)} \\ 
(14,1,2,11,28) & (14,1,8,5,28) & \underline{(2,1,8,21,32)}  & (8,1,2,21,32) & \underline{(8,1,16,7,32)}      \\    
(16,1,8,7,32) & (1,1,10,23,35) & (1,1,26,11,39) & \underline{(4,1,24,19,48)}  & (24,1,4,19,48)    \\
\underline{(12,1,16,19,48)} & (16,1,12,19,48) & (2,1,16,37,56) & (4,1,14,37,56) & (14,1,4,37,56)  \\
(2,1,40,17,60) & \underline{(12,1,30,17,60)} & (30,1,12,17,60) & (12,1,40,7,60) & (12,1,48,11,72) \\  
\underline{(24,1,36,11,72)} & (36,1,24,11,72) & (11,1,14,51,77) & \underline{(10,1,16,53,80)} & (16,1,10,53,80) \\
(3,1,54,23,81) & \underline{(8,1,44,35,88)}  & (44,1,8,35,88) & \underline{(8,1,26,69,104)}  & (26,1,8,69,104) \\
(15,1,70,19,105) & (6,1,96,41,144) & (7,1,46,107,161) & (84,1,16,67,168) & (50,1,16,133,200)  \\
(100,1,80,19,200) & (9,1,138,59,207) & (72,1,192,23,288)  & (14,1,88,205,308) & (18,1,264,113,396) \\
\hline
\end{tabular}
\label{tab:2c-1}
\end{table}
 
In the table of Egyptian fractions presented in Appendix~\ref{app:1}, the largest denominator
which occurs is 1806, and hence it makes sense to run a computer search using these bounds.  This gives 55 additional solutions, considering different possible covering maps (see table \eqref{tab:2c-1}).  There are 14 examples which have 2 different possible covering maps (we underline only the first).  So there are 41 underlying classes of hypersurfaces.

\subsection*{Case 2: $a=2$, $bq=2c$ and $2c-2=rd$.}

Here   \eqref{eq:fermat-equality} becomes
$$
    \frac{1}{2c} + \frac{1}{c} + \frac{1}{q} + \frac{d}{2c} + \frac{1}{2} = 1
$$
Moreover,
$ 
        2c-2 = rd $ implying $1 -1/c = rd/2c$, or, $d/2c = 1/r - 1/rc$
and hence
$$
\aligned
    1 &= \frac{1}{2c} + \frac{1}{c} + \frac{1}{q} + \frac{1}{r} - \frac{1}{rc} 
        +\frac{1}{2}
      = \frac{1}{2c} + \frac{1}{c}\left(\frac{r-1}{r}\right) + \frac{1}{q} + \frac{1}{r} + \frac{1}{2} \\ 
      &< \frac{1}{2c} + \frac{1}{c} + \frac{1}{q} + \frac{1}{r} + \frac{1}{2} 
\endaligned
$$
This produces another 5 examples, which are listed in Table~\ref{tab:2c-2}.
\begin{table}[h]\caption{$rd=2c-2$}
\begin{tabular}{|l|l|l|l|l|}
\hline
(1,1,2,3,7) & (1,1,2,6,10) & (5,1,2,7,15) & (13,1,2,10,26) & (7,1,2,18,28) \\
\hline
\end{tabular}
\label{tab:2c-2}
\end{table}

\subsection*{Other examples}
Searching a database of Calabi-Yau hypersurfaces yields
an additional 72 symmetric examples, which are listed in table \eqref{tab:misc-ex} (we did not try to enumerate all of the possible covers).

\begin{table}[h]\caption{Additional examples from a CY database} 
\begin{tabular}{|l|l|l|l|l|}
\hline
(1,1,3,5,10)       &  (1,1,3,7,12)     & (1,1,3,9,14)      & (1,1,4,8,14)      & (1,1,6,8,16)      \\
(1,1,6,10,18)      & (1,1,8,12,22)     & (1,1,8,20,30)     & (1,1,9,21,32)     & (1,1,11,15,28)    \\
(1,1,11,26,39)     & (1,1,12,16,30)    & (1,1,12,28,42)    & (2,1,4,5,12)      & (2,1,9,12,24)     \\
(3,1,3,11,18)      & (3,1,6,14,24)     & (9,1,3,23,36)     & (11,1,3,29,44)    & (3,1,12,20,36)    \\
(3,1,16,40,60)     & (3,1,24,32,60)    & (4,1,4,7,16)      & (8,1,4,11,24)     & (13,1,4,8,26)     \\
(11,1,4,28,44)     & (4,1,4,15,20,40)  & (5,1,6,18,30)     & (5,1,9,15,30)     & (18,1,5,12,36)    \\
(25,1,5,19,50)     & (29,1,5,23,58)    & (5,1,24,60,90)    &( 5,1,36,48,90)    & (1,6,7,14,28)     \\
(1,6,11,15,33)     & (18,1,6,11,36)    & (15,1,15,38,60)   & (7,1,8,12,28)     & (7,1,9,11,28)     \\
(7,1,16,32,56)     &  (9,1,8,36,54)    & (19,1,8,48,76)    & (8,1,27,36,72)    & (37,1,8,28,74)    \\
(45,1,8,36,90)     & (9,1,12,32,54)    & (20,1,9,30,60)    & (9,1,30,50,90)    & (33,1,10,22,66)   \\
(45,1,10,34,90)    & (35,1,11,93,140)  & (11,1,47,117,176) & (11,1,48,72,132)  & (20,1,12,27,60)   \\
(33,1,12,20,66)    & (27,1,12,68,108)  & (69,1,13,55,138)  & (13,1,83,111,208) & (15,1,64,160,240) \\ 
(15,1,96,128,240)  & (35,1,16,88,140)  & (51,1,16,136,204) & (69,1,16,52,138)  & (85,1,16,68,170)  \\
(17,1,24,60,102)   & (19,1,24,32,76)   & (81,1,19,61,162)  &  (53,1,20,32,106) & (85,1,20,64,170)  \\
(49,1,23,123,196)  & (51,1,24,128,204) &                   &                   &                   \\
\hline
\end{tabular}
\label{tab:misc-ex}
\end{table}

\vspace{22ex}
  \begin{appendix}
  \section{On Egyptian fractions}\label{app:1}
  \par In the tables below, we list all of the solutions to the equation
$$
    1=  \frac{1}{n} + \frac{1}{p} + \frac{1}{q} + \frac{1}{r} + \frac{1}{s}
$$
where $(n,p,q,r,s)$ are positive integers with $n\leq p\leq q\leq r\leq s$.  The only possible values of $n$ are $n=2,3,4,5$.   For $n=2$, we find 108 solutions.  For
$n=3$, we find 33 solutions.  For $n=4$, we find 5 solutions.  Finally, for 
$n=5$, the only solution is $1/5$ five times.  This gives 147 solutions
in total, which matches the value in the online encyclopedia of integers (A002966).
For each value of $n$, we list only possible values of $(p,q,r,s)$.

$$
\begin{array}{|l|l|l|l|l|l|}
\hline
{(3,7,43,1806)} & 
{(3,7,44,924)}  & 
{(3,7,45,630)}  & 
{(3,7,46,483)}  &
{(3,7,48,336)}  &
{(3,7,49,294)}  \\
{(3,7,51,238)}  & 
{(3,7,54,189)}  &
{(3,7,56,168)}  & 
{(3,7,60,140)}  & 
{(3,7,63,126)}  & 
{(3,7,70,105)}  \\
{(3,7,78,91)}   & 
{(3,7,84,84)}   &
{(3,8,25,600)}  & 
{(3,8,26,312)}  & 
{(3,8,27,216)}  & 
{(3,8,28,168)}  \\ 
{(3,8,30,120)}  & 
{(3,8,32,96)}   & 
{(3,8,33,88)}   &
{(3,8,36,72)}   &
{(3,8,40,60)}   & 
{(3,8,42,56)}   \\
{(3,8,48,48)}   & 
{(3,9,19,342)}  & 
{(3,9,20,180)}  &
{(3,9,21,126)}  &
{(3,9,22,99)}   & 
{(3,9,24,72)}   \\
{(3,9,27,54)}   & 
{(3,9,30,45)}   & 
{(3,9,36,36)}   &
{(3,10,16,240)} & 
{(3,10,18,90)}  & 
{(3,10,20,60)}  \\
{(3,10,24,40)}  & 
{(3,10,30,30)}  & 
{(3,11,14,231)} & 
{(3,11,15,110)} & 
{(3,11,22,33)}  & 
{(3,12,13,156)} \\
{(3,12,14,84)}  & 
{(3,12,15,60)}  & 
{(3,12,16,48)}  & 
{(3,12,18,36)}  & 
{(3,12,20,30)}  & 
{(3,12,21,28)}  \\ 
{(3,12,24,24)}  & 
{(3,13,13,78)}  & 
{(3,14,14,42)}  & 
{(3,14,15,35)}  & 
{(3,14,21,21)}  & 
{(3,15,15,30)}  \\
{(3,15,20,20)}  & 
{(3,16,16,24)}  & 
{(3,18,18,18)}  & 
{(4,5,21,420)}  & 
{(4,5,22,220)}  & 
{(4,5,24,120)}  \\
{(4,5,25,100)}  & 
{(4,5,28,70)}   & 
{(4,5,30,60)}   &
{(4,5,36,45)}   & 
{(4,5,40,40)}   & 
{(4,6,13,156)}  \\
{(4,6,14,84)}   & 
{(4,6,15,60)}   & 
{(4,6,16,48)}   & 
{(4,6,18,36)}   & 
{(4,6,20,30)}   & 
{(4,6,21,28)}   \\
{(4,6,24,24)}   & 
{(4,7,10,140)}  & 
{(4,7,12,42)}   & 
{(4,7,14,28)}   & 
{(4,8,9,72)}    &
{(4,8,10,40)}   \\ 
{(4,8,12,24)}   & 
{(4,8,16,16)}   & 
{(4,9,9,36)}    & 
{(4,9,12,18)}   & 
{(4,10,10,20)}  & 
{(4,10,12,15)}  \\
{(4,12,12,12)}  & 
{(5,5,11,110)}  & 
{(5,5,12,60)}   & 
{(5,5,14,35)}   & 
{(5,5,15,30)}   & 
{(5,5,20,20)}   \\  
{(5,6,8,120)}   & 
{(5,6,9,45)}    & 
{(5,6,10,30)}   & 
{(5,6,12,20)}   & 
{(5,6,15,15)}   & 
{(5,7,7,70)}    \\ 
{(5,8,8,20)}    & 
{(5,10,10,10)}  & 
{(6,6,7,42)}    & 
{(6,6,8,24)}    & 
{(6,6,9,18)}    & 
{(6,6,10,15)}   \\ 
{(6,6,12,12)}   & 
{(6,7,7,21)}    & 
{(6,8,8,12)}    &
{(6,9,9,9)}     & 
{(7,7,7,14)}    &
{(8,8,8,8)}\\
\hline 
\end{array}
$$
\centerline{\textbf{First table  $n=2$}  i.e., $
     1 = \frac{1}{2} + \frac{1}{p} + \frac{1}{q} + \frac{1}{r} + \frac{1}{s}
 $}
%

$$
\begin{array}{|l|l|l|l|l|l|}
\hline
(3,4,13,156) & 
(3,4,14,84)  & 
(3,4,15,60)  & 
(3,4,16,48)  & 
(3,4,18,36)  & 
(3,4,20,30) \\
(3,4,21,28) & 
(3,4,24,24) &
(3,5,8,120) & 
(3,5,9,45)  & 
(3,5,10,30) & 
(3,5,12,20) \\
(3,5,15,15) & 
(3,6,7,42)  & 
(3,6,8,24)  & 
(3,6,9,18)  & 
(3,6,10,15) & 
(3,6,12,12) \\
(3,7,7,21)  &
(3,8,8,12)  & 
(3,9,9,9)   & 
(4,4,7,42)  & 
(4,4,8,24)  & 
(4,4,9,18)  \\
(4,4,10,15) & 
(4,4,12,12) &
(4,5,5,60)  & 
(4,5,6,20)  & 
(4,6,6,12)  & 
(4,6,8,8)   \\
(5,5,5,15)  & 
(5,5,6,10)  & 
(6,6,6,6)   & 
& & \\
\hline
\end{array}
$$
\centerline{\textbf{Second table  $n=3$}  i.e., $
     1 = \frac{1}{3} + \frac{1}{p} + \frac{1}{q} + \frac{1}{r} + \frac{1}{s}
      $}

$$
\begin{array}{|l|l|l|l|l|}
\hline
(4,4,5,20) & 
(4,4,6,12) & 
(4,4,8,8)  & 
(4,5,5,10) & 
(4,6,6,6)  \\
\hline
\end{array}
$$
\centerline{\textbf{Third table  $n=4$}  i.e., $
     1 = \frac{1}{4} + \frac{1}{p} + \frac{1}{q} + \frac{1}{r} + \frac{1}{s}
      $}

\section{SAGE CODE for Table~\ref{tab:First}} \label{app:2}

\begin{sageblock}
% Code to create the table; one needs to load J. Carlson's Sage code from
% https://github.com/jxxcarlson/math_research/blob/master/hodge.sage
%
A=3; a=4; b=16; c=24; # Need 2c divisible by A.
print(A,1,a,b,c)
m=2*c/A;
for d in range(2,m+1):
    if ((m % d)==0):
       s=0; q=0;
       for k in range(1,d):
           if gcd(k,d)==1:
      #hodge(2*c,[1,a,b,c],m,k*(m/d))
s = s + (hodge(2*c,[1,a,b,c],m,k*(m/d)))[1]
q = q + (hodge(2*c,[1,a,b,c],m,k*(m/d)))[2]
       if(d<m):
print((s+q)/euler_phi(d),d)
       if(d==m):
print((2+s+q)/euler_phi(d),d)
\end{sageblock}
%

  \end{appendix}

 \bibliographystyle{plain}
\bibliography{/Users/chrismacbook/surfdriveLeiden/Werk/Bib.f/bisect.bib} 

\end{document}